\documentclass[11pt]{article}
\usepackage[T1]{fontenc}
\usepackage[utf8]{inputenc}
\usepackage{lmodern}
\usepackage[margin=1in]{geometry}
\usepackage{amsmath,amssymb}
\usepackage{graphicx}
\usepackage{booktabs}
\usepackage{xcolor}
\usepackage{url}
\usepackage[hidelinks]{hyperref}
\usepackage{natbib}
\usepackage[section]{placeins}

\definecolor{wsp}{HTML}{2A78D6}
\newcommand{\finding}[2]{\par\smallskip\noindent\textbf{#1.}~#2\par\smallskip}

\title{Can a Language Model Write the Weights?\\
\large An Empirical Study of LLM-Written Scoring Networks for Column Generation and MILP Branching,\\
with a Bus-Driver Scheduling Application}

\author{Renan Butkeraites\\
BIS --- Butkeraites Intelligent Solutions\\
\texttt{renan@butkeraites.com}}

\date{November 2026}

\begin{document}
\maketitle

\begin{abstract}
\emph{Weight-space programming} (WSP) is the idea that a large language model (LLM) can write the
weights of a small neural network directly, guided only by reading the outputs of forward passes,
with no gradient training \citep{wsp}. We test whether WSP can improve classical, exact operations
research algorithms in places where they rely on cheap scoring rules. In four studies, an LLM wrote
small multilayer perceptrons that (i) select arcs, labels and ng-neighbourhoods inside column
generation for the vehicle routing problem with time windows, and (ii) choose branching variables
in SCIP's branch-and-bound. The algorithms stay exact whatever the weights are; a bad network only
costs time. Results: in column generation, no lever gave a material gain (exact pricing dominates
the running time). In branching on synthetic benchmark classes, a hand-written network solved
independent-set instances in 0.37$\times$ the time of SCIP's default rule (18 paired instances,
Wilcoxon $p<10^{-5}$) and matched or beat networks trained by imitation of strong branching on the
same features. The same rule was not general: on 39 MIPLIB 2017 instances it solved 19 versus 24.
In an application built from the real GTFS feed of Atlanta's transit agency (MARTA), a crew
scheduling set-partitioning model solved by a native C plugin with SCIP's factory settings, the
LLM-written rules took 0.63$\times$ the default's time on 32 held-out instances ($p<10^{-4}$) ---
but SCIP's own plain pseudocost rule did exactly as well (ratio 0.99, $p=0.19$). We conclude that,
in these settings, the gains came from being cheap, not from what the LLM wrote. All code, logs
and the numbers behind every figure are released.
\end{abstract}

\begin{center}
\fbox{\begin{minipage}{0.94\linewidth}\small
\textbf{AI-generation and provenance statement (read first).}
This paper --- its text, figures, tables, code, experimental design and analysis --- was generated
by an AI model, \textbf{Claude Opus~5.5 (Anthropic)}, running as an agent in Claude Code, under the
direction of the human author. The human author chose the research questions, approved or
redirected each step, and takes full responsibility for the content; the AI is not an author.
\textbf{All numerical results are from computations actually run} on the author's local machine
(Apple M5, 10 cores, 32\,GB RAM, macOS~27) between 2026-11-03 and 2026-11-05; no number was
estimated, simulated or typed in by the model. Every number in the paper is regenerated from the
raw run logs by one script (\texttt{scripts/figs.py}, Appendix~\ref{app:repro}). Every reference
was checked against the publisher's page, a DOI resolver, dblp or the original web page.
Note the circularity: the same model that wrote this paper also wrote the network weights that are
the object of study. Appendix~\ref{app:ai} lists what the model did and what the human did.
\end{minipage}}
\end{center}

\section{Introduction}

Exact algorithms in operations research (OR) contain many small decisions that are made by cheap
scoring rules: which arcs to keep in a pricing graph, which labels to extend, which fractional
variable to branch on. These rules are where machine learning has entered the field
\citep{bengio2021}: branching rules have been learned by ranking \citep{khalil2016}, by
regression on strong-branching scores \citep{alvarez2017} and by graph neural networks trained to
imitate strong branching \citep{gasse2019}. All of them need training data and a training
pipeline.

A recent post, ``Sparks of Weight-Space Programming'' \citep{wsp}, reports that an LLM can instead
\emph{write} the weights of a small network by hand, reading the outputs of forward passes and
editing the weights, without any gradient step. If that works for OR scoring rules, it would give
a third option between hand-coded heuristics and trained models: a network whose every weight has
a written justification, specialised to a problem family in hours, with no training set.

This paper asks a narrow, testable question:

\begin{quote}
\emph{Can an LLM write the weights of a cheap scoring network that makes a gold-standard exact OR
algorithm faster, on held-out instances, without changing its result?}
\end{quote}

\begin{figure}[!htbp]
\centering
\includegraphics[width=0.92\linewidth]{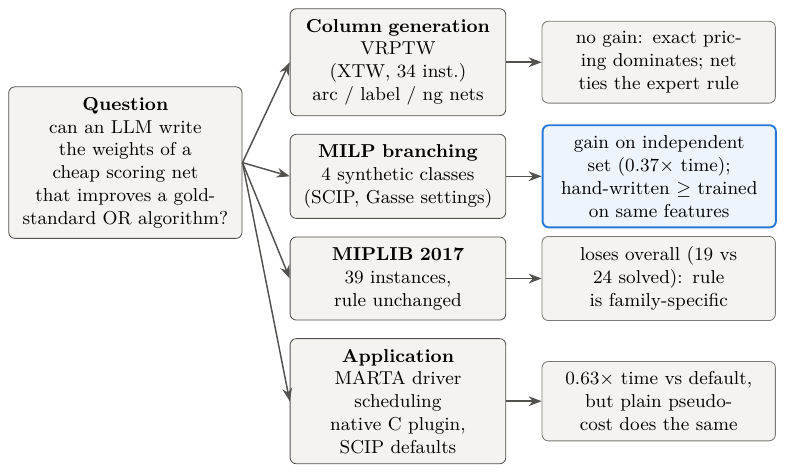}
\caption{Map of the study. Each row is one study; the right column is its verdict. Only the
synthetic branching study (blue) found a gain attributable to the written weights, and it did not
transfer to MIPLIB or to the application.}
\label{fig:map}
\end{figure}

We ran four studies (Figure~\ref{fig:map}). The setting is deliberately favourable to measurement:
the host algorithms are exact, so the network can only change time and tree size, never the
optimum; and we always compare against the incumbent expert rule in the same batch of runs.

\paragraph{Findings at a glance.} Each finding has an identifier used in the claims table
(Appendix~\ref{app:claims}), which points to the figure, table and log file that support it.
\finding{F1}{Column generation for VRPTW: no material gain. Exact pricing takes most of the time,
and the LLM-written arc scorer ties the classical reduced-cost rule (time ratio
1.02 over 34 instances; Section~\ref{sec:cg}).}
\finding{F2}{Synthetic MILP branching: on independent set, the LLM-written rule needs
0.37$\times$ the time of SCIP's default with similar or smaller trees ($n=18$, $p<10^{-5}$); on
set cover it ties (0.90$\times$, $p=0.14$) with larger trees (Section~\ref{sec:syn}).}
\finding{F3}{With the same 7 features, hand-written weights match or beat weights trained by
imitation of strong branching in 3 of 4 settings, although the trained weights imitate strong
branching better (Section~\ref{sec:trained}).}
\finding{F4}{The rule is family-specific: on 39 MIPLIB 2017 instances it solves 19 vs 24
(Section~\ref{sec:miplib}).}
\finding{F5}{Application (MARTA bus-driver scheduling, native C plugin, SCIP defaults, 32 held-out
instances): LLM-written rules take 0.63$\times$ the default's time, but SCIP's plain pseudocost
rule achieves the same (0.99$\times$ relative to it, $p=0.19$). The gain is ``skip strong
branching'', not WSP (Section~\ref{sec:marta}).}
\finding{F6}{A net-guided lookahead (strong branching only on the network's top-$k$) halved tree
sizes relative to SCIP's default on 8 practice instances, but did not generalise to the test set
(1.05$\times$ time, $p=0.51$).}

\section{Background}

\paragraph{Column generation.} Column generation solves a linear relaxation with exponentially
many columns by alternating a restricted master LP and a pricing problem that finds columns with
negative reduced cost \citep{lubbecke2005,desaulniers2005}. For the vehicle routing problem with
time windows (VRPTW) the pricing problem is an elementary shortest path problem with resource
constraints \citep{desrochers1992}, commonly relaxed to ng-routes \citep{baldacci2011}. Practical
codes price heuristically on a reduced graph first and fall back to exact pricing; dual smoothing
\citep{wentges1997} stabilises the duals.

\paragraph{Branching.} In LP-based branch-and-bound, the choice of the branching variable drives
tree size. Full strong branching solves the two child LPs for every candidate and picks the best
product of gains; it gives small trees at a high cost per node. Reliability pseudocost branching,
SCIP's default, uses strong branching only until pseudocost estimates are reliable
\citep{achterberg2005,achterberg2007}. Learning to branch tries to approximate strong branching
cheaply \citep{khalil2016,alvarez2017,gasse2019}. We use SCIP \citep{achterberg2009,scip10} through
PySCIPOpt \citep{pyscipopt} and natively in C.

\begin{figure}[!htbp]
\centering
\includegraphics[width=0.72\linewidth]{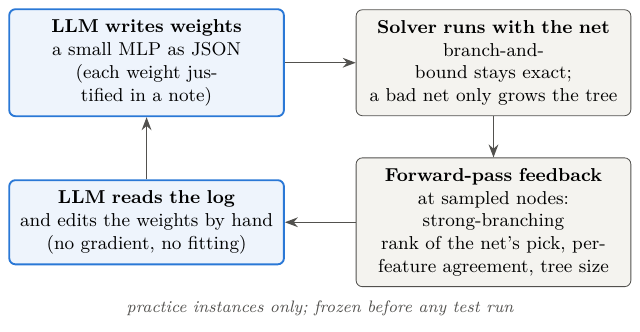}
\caption{The WSP protocol used in all studies. Weights are written by the LLM as a JSON file
with a free-text note justifying each weight; feedback runs are read, never fitted; the network is
frozen before the held-out test.}
\label{fig:loop}
\end{figure}

\paragraph{Weight-space programming.} In WSP \citep{wsp} the LLM is the optimiser: it proposes
weights, observes forward-pass outputs, and edits the weights in text. Our protocol (Figure~\ref{fig:loop})
follows that spirit: diagnostics are \emph{read} by the model and never fitted.

\section{Method}

\begin{figure}[!htbp]
\centering
\includegraphics[width=0.95\linewidth]{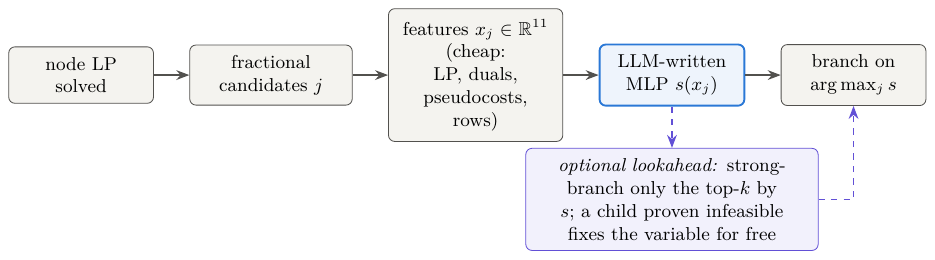}
\caption{The branching rule at one node. The dashed path is the optional net-guided lookahead
(Section~\ref{sec:marta}).}
\label{fig:bnb}
\end{figure}

\paragraph{Where the network sits.} In branching (Figure~\ref{fig:bnb}), at every node SCIP
calls our branching rule with the fractional candidates. We compute a feature vector per
candidate from quantities already available at the node (Table~\ref{tab:features}), evaluate the
network, and branch on the arg max. The network is a multilayer perceptron stored as JSON; most
networks in this paper have one linear layer or one ReLU layer of width 6. An optional
\emph{lookahead} runs strong branching only on the network's top-$k$ candidates, and turns any
child proven infeasible into a free bound change, as SCIP's own strong branching does.

\begin{table}[!htbp]
\centering\small
\caption{Candidate features. Features marked $^\dagger$ are divided by their maximum over the
node's candidates. The last four were added for set-partitioning structure in the application.}
\label{tab:features}
\begin{tabular}{@{}lp{0.72\linewidth}@{}}
\toprule
Feature & Definition (candidate $j$ with LP value $v_j$, fractional part $f_j$) \\
\midrule
\texttt{frac} & $2\min(f_j, 1-f_j)$; 1 means exactly one half \\
\texttt{ps}$^\dagger$ & SCIP pseudocost score of $j$ \\
\texttt{obj}$^\dagger$ & $|c_j|$ \\
\texttt{deg}$^\dagger$ & number of LP rows in which column $j$ has a nonzero \\
\texttt{tight} & share of $j$'s rows that are tight at the LP optimum \\
\texttt{dualm}$^\dagger$ & $\sum_i |y_i a_{ij}|$ with $y$ the row duals \\
\texttt{val} & $v_j$ \\
\texttt{fconf}$^\dagger$ & LP mass of the other fractional candidates sharing a row with $j$ \\
\texttt{alt} & $1/(m-1)$, $m$ the fewest fractional candidates in any row of $j$ \\
\texttt{rowf}$^\dagger$ & number of $j$'s rows touched by more than one fractional candidate \\
\texttt{eff} & cost per row of the most efficient candidate divided by that of $j$ \\
\bottomrule
\end{tabular}
\end{table}

\paragraph{Protocol rules.} (1) Separate practice and test instance sets generated with different
seeds; weights are written on practice only and frozen before the test. (2) Every comparison is
\emph{paired}: all rules run on the same instances in the same batch, so they share machine load.
(3) Tree size (nodes) is deterministic for a given SCIP build and seed; time is noisy across
batches, so time is only compared within a batch. (4) Summary statistics are shifted geometric
means (shift 10 for nodes, 1\,s for time) with time-outs counted at the limit. (5) Paired
two-sided Wilcoxon signed-rank tests on log time ratios; pairs that both hit the time limit are
excluded from the test.

\paragraph{Software.} SCIP~10.0 through PySCIPOpt~6.3.0 for Studies 2--3 and the Python
prototype of Study~4; SCIP~10.1.0 (SoPlex~8.1.0) linked from C for the final application runs;
HiGHS \citep{huangfu2018} through the \texttt{highs} Rust crate~2.4.0 for column generation;
Python~3.12, Rust~1.96. Studies~2--3 use the settings of \citet{gasse2019} (no restarts, cuts only
at the root, fixed seed); the final application runs use SCIP's factory settings.

\section{Study 1: column generation for VRPTW}
\label{sec:cg}

\paragraph{Setup.} We implemented column generation for VRPTW in Rust: a set-covering master LP
kept warm in HiGHS, ng-route labelling with ng-sets of size 8, Wentges smoothing, heuristic
pricing on a reduced graph with a label cap, and exact pricing as fallback, so the final bound is
exact for every strategy. Instances are the XTW set released on CVRPLIB \citep{xtw}; 34 instances
with up to 256 customers for the main comparison and the 10 largest (819--1001 customers) for a
scale test. The LLM wrote networks for three levers: (a) arc scoring for the reduced graph,
compared with the classical rule ``keep the $K$ arcs with lowest reduced cost per customer'';
(b) label eviction under the cap; (c) ng-neighbourhood selection.

\begin{figure}[!htbp]
\centering
\includegraphics[width=0.92\linewidth]{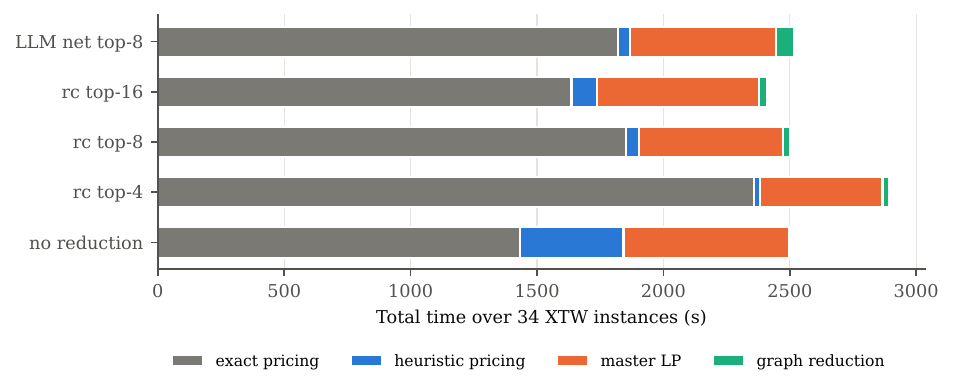}
\caption{Where column-generation time goes, summed over 34 XTW instances, for each arc-reduction
strategy. Exact pricing dominates in every strategy, so better heuristic arcs can save little; the
LLM-written net (top) ties the reduced-cost top-8 rule.}
\label{fig:cg}
\end{figure}

\paragraph{Results.} Figure~\ref{fig:cg} shows the time breakdown. The arc network ties the
reduced-cost top-8 rule (geometric mean time ratio 1.02). Label-eviction networks did not beat
``keep the cheapest labels'', and raising the label cap from 8 to 64 cut exact-pricing time by
only about 13\% while multiplying heuristic time about tenfold, so the eviction policy had little
headroom. ng-set selection that drops neighbours unable to close a cycle raised the bound by
0.03--0.15\% on average; a replicated timing at ng-8 gave 0.93$\times$ time, within noise. On the
largest instances, in heuristic-only mode, the arc network reached a lower master objective than
the top-8 rule at equal iteration on 7 of 10 instances, but its unoptimised scoring cost
40--197\,s, so it lost at equal time. \textbf{Verdict (F1):} no material gain.

\section{Study 2: branching on synthetic classes}
\label{sec:syn}

\paragraph{Instances.} We generated the four classes used by \citet{gasse2019}: set cover
\citep{balas1980} (1000 rows $\times$ 1000 columns, density 0.05), combinatorial auction
\citep{leytonbrown2000} (100 items, 1500 bids), capacitated facility location
\citep{cornuejols1991} (70 customers, 70 facilities) and maximum independent set (750 nodes, affinity 4) with a clique-partition
formulation, following the procedure that \citet{gasse2019} attribute to \citet{bergman2016};
like their public code, we draw the graphs from the Barabási--Albert model \citep{barabasi1999}. Practice: 4 per class; test: 10 per class; medium (larger, new
seed): 8 per class for set cover (1250 rows) and independent set (1000 nodes); a separate training
set of 10 per class for the trained baseline. Time limit 300\,s.

\begin{figure}[!htbp]
\centering
\includegraphics[width=0.95\linewidth]{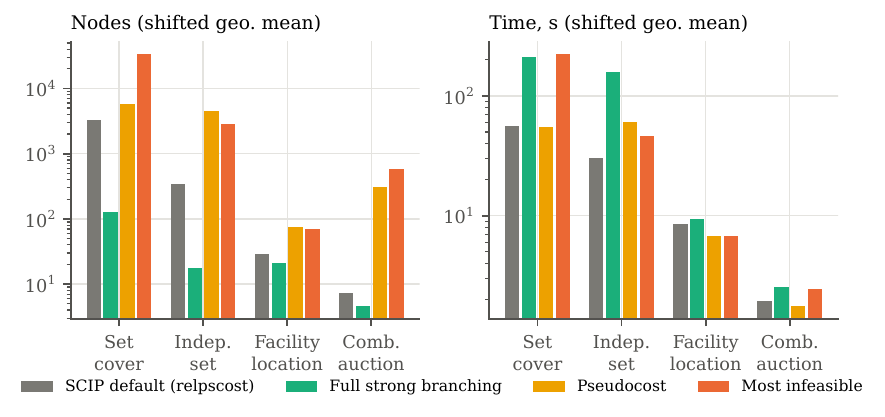}
\caption{Headroom on the synthetic test set (10 instances per class, 300\,s). Full strong
branching builds far smaller trees on set cover and independent set, at a much higher cost per
node.}
\label{fig:headroom}
\end{figure}

\paragraph{Headroom first.} Before writing any weights we measured the gap between SCIP's default
and full strong branching (Figure~\ref{fig:headroom}). On set cover and independent set, strong
branching trees are 26$\times$ and 19$\times$ smaller than the default's but about 4--5$\times$ slower,
which is the regime where a cheap, well-aimed rule can pay off. Facility location and auctions
have trees of tens of nodes and little headroom.

\paragraph{The written networks.} The LLM read per-feature agreement with strong branching on the
practice set and wrote network v0, then edited it into v1 after reading where v0's picks disagreed
with strong branching. In closed form, v1 scores a candidate as
\begin{equation}
s_{\text{v1}}(x) = \texttt{deg} + 0.7\,\texttt{dualm} + 0.4\,\texttt{val} + 0.2\,\texttt{ps} +
0.3\,\texttt{frac} - 3\,\max(0,\ 0.4 - \texttt{frac}).
\label{eq:v1}
\end{equation}
The note stored with the weights explains each term: \texttt{deg} favours hub variables whose
up-branch propagates widely; \texttt{dualm} measures how much LP dual pressure the column carries;
the ReLU gate removes nearly integral candidates; the pseudocost weight was cut from 0.7 to 0.2
because the disagreements came from candidates with high but noisy early pseudocosts.

\begin{figure}[!htbp]
\centering
\includegraphics[width=0.95\linewidth]{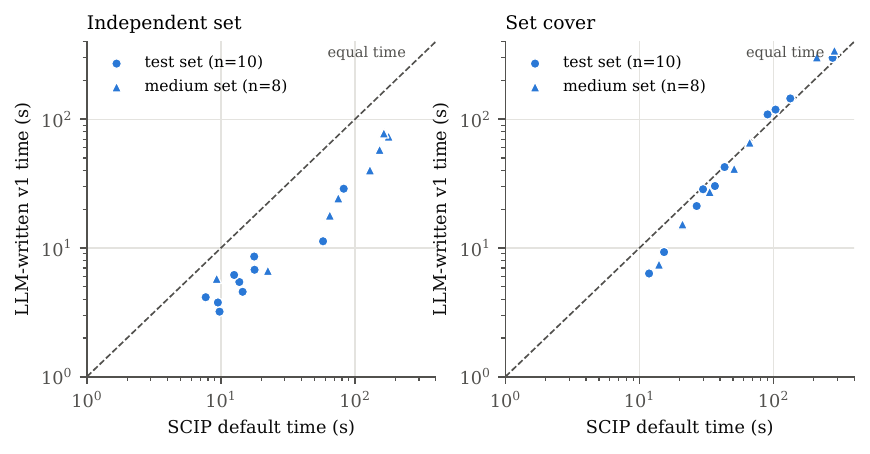}
\caption{Per-instance time, LLM-written v1 versus SCIP default, paired runs, test and medium sets.
Points below the diagonal are instances where v1 is faster. Independent set: 18 of 18 faster.}
\label{fig:syn}
\end{figure}

\paragraph{Results.} Figure~\ref{fig:syn}. On independent set, v1 is faster on all 18 instances,
with a geometric-mean time ratio of 0.37 ($p=7.6\times10^{-6}$) and smaller trees on the medium
set (1426 vs 2786 nodes). On set cover it ties on time (0.90, $p=0.14$) with 1.3--1.4$\times$ larger
trees. The network's features are computed in Python while SCIP's rules run in C, so these time
ratios are pessimistic for the network. Optimal values were identical on every instance solved by
both rules. \textbf{(F2)}

\subsection{Hand-written versus trained weights}
\label{sec:trained}

\begin{figure}[!htbp]
\centering
\includegraphics[width=0.92\linewidth]{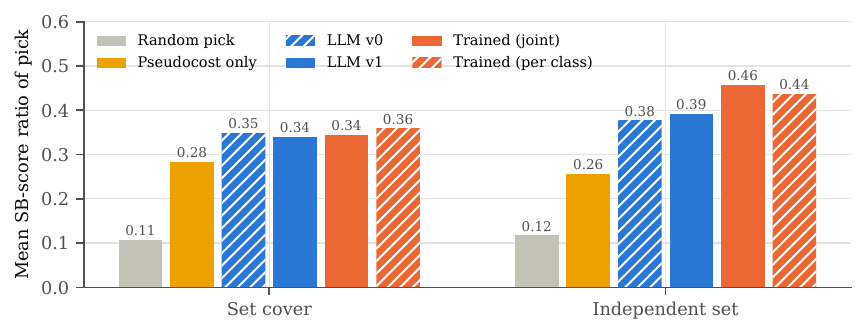}
\caption{Imitation quality on the training nodes: mean ratio between the strong-branching score
of the picked candidate and the best one (1 = always picks the strong-branching best). Trained
networks imitate strong branching better on independent set.}
\label{fig:agree}
\end{figure}

\begin{figure}[!htbp]
\centering
\includegraphics[width=0.95\linewidth]{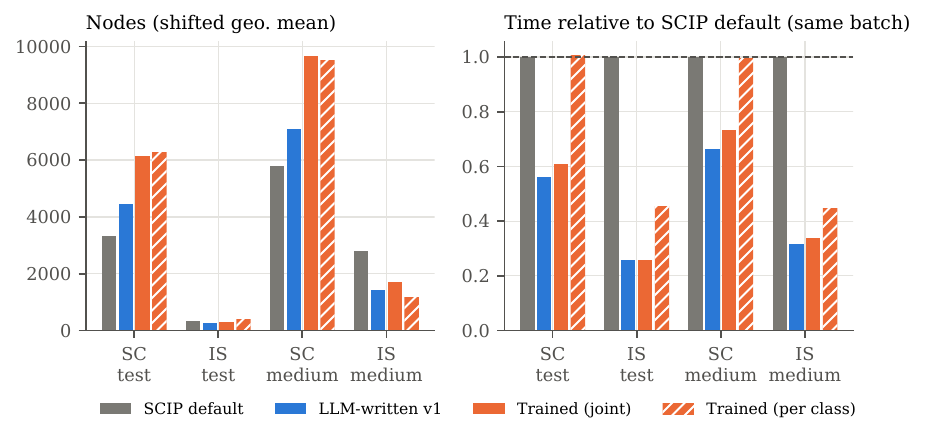}
\caption{Hand-written v1 versus trained networks with identical features, paired batch. Left:
tree size. Right: time relative to SCIP's default in the same batch. SC = set cover, IS =
independent set.}
\label{fig:trained}
\end{figure}

To separate ``LLM-written'' from ``a network on these features'', we trained the same
architecture (7 inputs, 16 ReLU units) by imitation of strong branching: listwise softmax
cross-entropy with the strong-branching best candidate as target, Adam, 20\% of instances held
out, on 2,166 labelled nodes (1,361 set cover, 805 independent set) from the separate training
set. Figure~\ref{fig:agree} shows imitation quality; Figure~\ref{fig:trained} shows the paired
solving results.

The trained networks imitate strong branching better on independent set (0.44--0.46 vs 0.39),
yet produce larger trees in 3 of 4 settings and are never faster than v1. Better imitation did not
mean better decisions. \textbf{(F3)} This is a deliberately small baseline (same 7 features, 2k
nodes), not a reproduction of the graph network of \citet{gasse2019}.

\section{Study 3: generalisation to MIPLIB 2017}
\label{sec:miplib}

\begin{figure}[!htbp]
\centering
\includegraphics[width=0.92\linewidth]{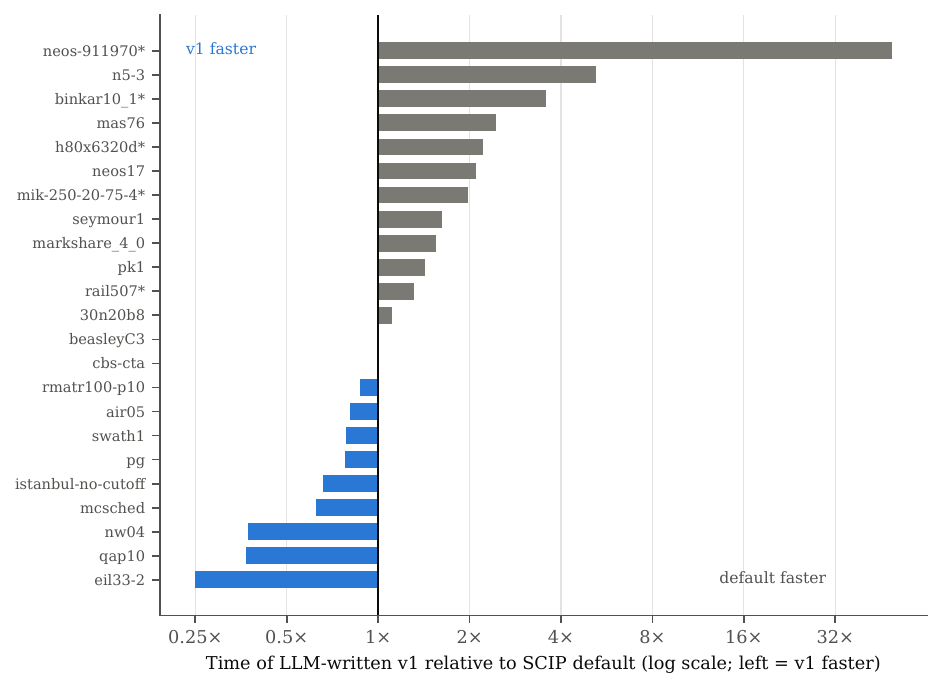}
\caption{MIPLIB 2017: time of v1 relative to SCIP's default, per instance (instances where both
rules hit the time limit, or both finish within 10\,s, are omitted). $^*$ = at least one rule hit
the 300\,s limit.}
\label{fig:miplib}
\end{figure}

We ran v1, unchanged, on 40 instances of MIPLIB 2017 \citep{miplib2017} chosen from the benchmark
set to be solvable within minutes (one could not be read), paired with SCIP's default, 300\,s,
same settings. Figure~\ref{fig:miplib} shows per-instance time ratios.

v1 solves 19 instances against 24 for the default; on the 17 non-trivial instances solved by both,
its trees are 2.6$\times$ larger (5080 vs 1989 nodes). It wins clearly on airline crew scheduling
(\texttt{nw04}), a quadratic assignment instance (\texttt{qap10}) and a routing instance
(\texttt{eil33-2}), and loses badly on others (\texttt{neos-911970}: 6\,s for the default, time-out
for v1). A rule written for two synthetic families is not a general-purpose rule. \textbf{(F4)}

\section{Study 4: application to bus-driver scheduling (MARTA)}
\label{sec:marta}

The pitch for WSP in practice is per-family specialisation: a company that solves the same model
every day could have a rule written for its family. We tested this on driver scheduling
(run-cutting), a classic transit problem \citep{wren1995} usually written as set partitioning
\citep{ryan1981}.

\begin{figure}[!htbp]
\centering
\includegraphics[width=\linewidth]{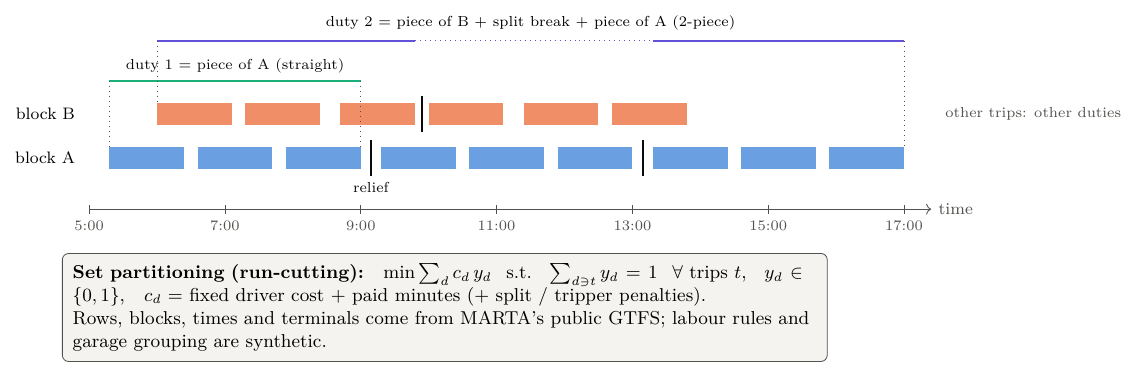}

\medskip
\includegraphics[width=0.42\linewidth]{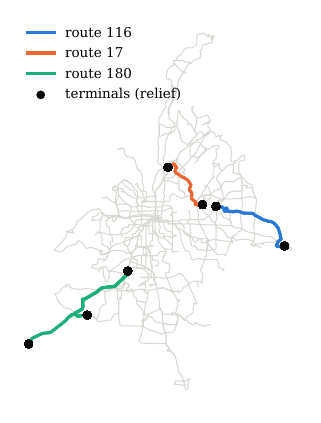}
\caption{Top: from vehicle blocks to duties to the set-partitioning model. Bottom: one practice
garage (three MARTA bus routes, coloured) over the whole bus network (grey); dots are terminals
where relief can happen.}
\label{fig:marta}
\end{figure}

\paragraph{Instances from real data.} From the public GTFS feed of the Metropolitan Atlanta Rapid
Transit Authority \citep{gtfs} we took every bus trip with its vehicle block, times and terminals
(6,918 weekday bus trips on 81 routes). A ``garage'' is a set of three routes whose blocks stay within them
(GTFS does not publish garage assignment). Within each block, a \emph{piece} is a run of
consecutive trips; relief (driver change) happens at trip boundaries. A duty is one piece of
5--8.5\,h of work (straight), two pieces with a 30--90\,min meal break or a 1.5--5\,h split break
(spread $\le$ 12\,h, relief within 400\,m), or a short ``tripper'' piece with a penalty
(Figure~\ref{fig:marta}). Costs: a fixed cost per driver, paid minutes with an 8\,h guarantee,
and split and tripper penalties. At most 15 second pieces are sampled per first piece and
instances above 12,000 duties are discarded. The labour rules and the garage grouping are
synthetic; trips, blocks, times and terminals are real. Practice: 8 non-trivial garages from a
pool of 40 (seed 7); test: 32 new garages (seed 11), weekday, alternate weekday, Saturday and
Sunday service, used without any filtering.

\begin{figure}[!htbp]
\centering
\includegraphics[width=0.92\linewidth]{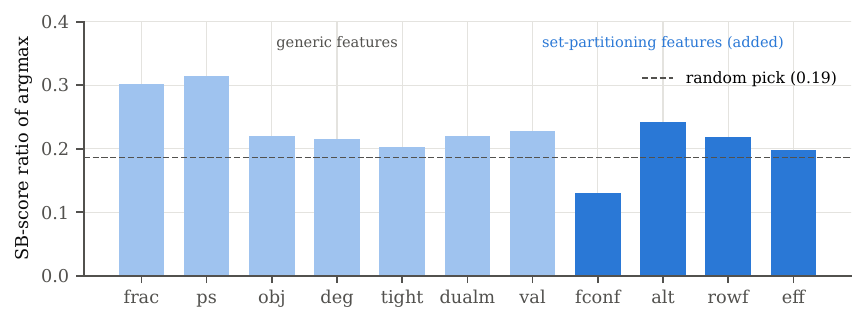}
\caption{MARTA practice diagnostics (2,264 nodes): strong-branching score ratio obtained by
branching on the arg max of each single feature. The dashed line is a random pick.}
\label{fig:mfeat}
\end{figure}

\paragraph{Headroom and diagnosis.} Strong branching solved practice garages with 11--75 nodes
against 286--2401 for the default (SCIP~10.0, Gasse settings). The per-feature diagnostics
(Figure~\ref{fig:mfeat}) showed that no generic feature is informative and that in 91\% of nodes
some candidate has a child that strong branching proves infeasible or worse than the incumbent:
the good decision is ``which variable can be fixed for free'', which cheap features cannot see.
The LLM added the four set-partitioning features of Table~\ref{tab:features}; they were no better.

\begin{figure}[!htbp]
\centering
\includegraphics[width=0.6\linewidth]{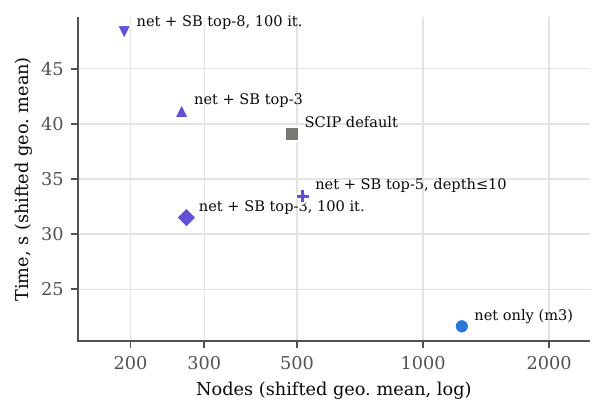}
\caption{MARTA practice set (8 garages, native C, SCIP factory settings): tree size versus time
for SCIP's default, the network alone, and network-guided lookahead variants. Two rules were
frozen from this plot: \texttt{vm\_fast} (network only) and \texttt{vm\_la} (top-3, 100
iterations).}
\label{fig:mprac}
\end{figure}

\paragraph{Forward passes and the lookahead.} The LLM wrote four linear candidates and measured
tree size and time on practice; the best (\texttt{m3}: $\texttt{frac} + 0.5\,\texttt{ps} +
0.3\,\texttt{val} + 0.3\,\texttt{alt} + 0.5\,\texttt{rowf}$) was then combined with a net-guided
lookahead (Figure~\ref{fig:bnb}). To make time comparable with SCIP's C rules, the rule was ported
to a native SCIP plugin in C (\texttt{wspsolve.c}); the C port reproduces the Python node counts
(1226 vs about 1290 on one garage). On practice, with factory settings (Figure~\ref{fig:mprac}),
strong branching on the network's top-3 with a 100-iteration limit produced trees 0.56$\times$ the
default's and ran in 0.81$\times$ its time; top-8 gave 0.40$\times$ the nodes.

\begin{figure}[!htbp]
\centering
\includegraphics[width=0.92\linewidth]{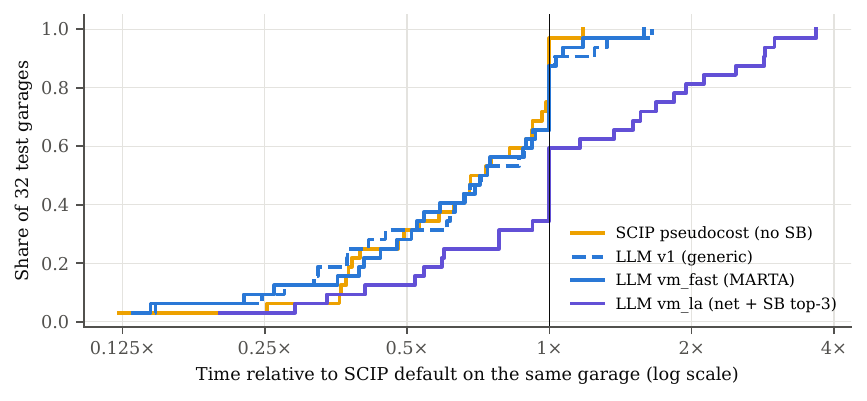}
\caption{MARTA held-out test (32 garages): cumulative distribution of the per-garage time ratio
against SCIP's default in the same batch. Curves to the left are faster. SCIP's plain pseudocost
rule and the LLM-written rules overlap.}
\label{fig:mtest}
\end{figure}

\paragraph{Held-out test.} Figure~\ref{fig:mtest} and Table~\ref{tab:marta}. Against SCIP's
default, the LLM-written rules without lookahead are faster on about 70\% of garages, with a
geometric-mean time ratio of 0.63 ($p=4\times10^{-5}$), and solve as many or more instances. The
lookahead variant did not generalise (1.05$\times$, $p=0.51$, 22 solved). The decisive control is
SCIP's own pseudocost rule, which uses no strong branching at all: it is as fast as the
LLM-written rules (\texttt{vm\_fast}/\texttt{pscost} time ratio 0.99, $p=0.19$), and the
MARTA-specific network is no better than the generic v1 (ratio 1.00, $p=0.61$). \textbf{(F5, F6)}

\begin{table}[!htbp]
\centering\small
\caption{MARTA held-out test, native C, SCIP~10.1 factory settings, 300\,s, paired. Shifted
geometric means over all 32 garages (time-outs at the limit). The first four rows share one batch;
the last two rows come from a second batch with its own default run.}
\label{tab:marta}
\begin{tabular}{@{}lrrr@{}}
\toprule
Rule & Solved & Nodes & Time (s) \\
\midrule
SCIP default (relpscost) & 24 & 798 & 88.9 \\
SCIP pseudocost & 24 & 2074 & 57.2 \\
LLM v1 (generic, from Study 2) & 25 & 1987 & 56.9 \\
LLM \texttt{vm\_fast} (MARTA) & 24 & 1980 & 56.6 \\
\midrule
SCIP default (second batch) & 23 & 782 & 88.3 \\
LLM \texttt{vm\_la} (net + top-3 lookahead) & 22 & 640 & 93.8 \\
\bottomrule
\end{tabular}
\end{table}

\section{Discussion}

\paragraph{What the written weights bought.} One clear, attributable gain: independent set in
Study~2, where v1 beat both SCIP's default and SCIP's pseudocost rule (Figure~\ref{fig:headroom}:
pseudocost is 2$\times$ slower than the default there, v1 is 2.7$\times$ faster; different batches), and beat trained
networks with the same inputs. Its terms are interpretable and were written from a few diagnostic
reads. Everywhere else, the gain over the default was explained by a cheaper rule that SCIP
already ships, or did not exist.

\paragraph{Why.} Strong-branching agreement measured from cheap features is low for every feature
and every network we tried (0.2--0.46); the structure that makes strong branching valuable in
set partitioning (children that can be pruned) is not visible to cheap features. The only route
that reached strong-branching-like trees was spending strong branching selectively, and its
selection, tuned on 8 practice instances, overfit.

\paragraph{Threats to validity.} (i) Runs were concurrent (up to 10 processes on 10 cores);
pairing controls load within a batch but absolute times vary across batches. (ii) Studies~2--3
compute features in Python. (iii) Two SCIP versions were used (10.0 via PySCIPOpt, 10.1 natively);
the native default solved one practice garage in 21 nodes versus 337 in 10.0, so conclusions are
drawn within one build. (iv) Instance sets are small (8--40 per comparison). (v) The MARTA labour
rules and garage grouping are synthetic. (vi) Practice instances were chosen as non-trivial for
the default. (vii) The trained baseline is small. (viii) The same AI wrote the weights, ran the
analysis and wrote this paper; to limit this, every number is produced by a script from logs that
are released, and the negative results are reported as found.

\section{Conclusion}

An LLM can write branching networks that are readable, exact-safe and occasionally better than a
mature solver's default on a specific family, and better than gradient-trained weights on the same
features. In our application, a real-data bus-driver scheduling model, the measured speed-up of
about 37\% was real but was entirely reproduced by switching SCIP to its plain pseudocost rule, a
one-line configuration change. For practitioners the lesson is to test the cheapest baseline that
the solver already offers before attributing a gain to a learned or written rule. For WSP, the
open question is whether richer, problem-aware inputs (or a written policy for \emph{when} to spend
strong branching) can give the network information that cheap rules lack.

\section*{Data and code availability}
All code (Rust column generation, Python and C branching rules, MARTA instance generator), network
weight files with their notes, raw run logs (JSONL) and the figure script are in the ancillary
files of this submission and at \url{https://github.com/butkeraites/llm-written-weights-or}. Third-party data: XTW instances
from CVRPLIB \citep{xtw}, MIPLIB 2017 \citep{miplib2017}, MARTA GTFS \citep{gtfs}.

\bibliographystyle{plainnat}
\bibliography{refs}

\appendix

\section{Reproducibility}
\label{app:repro}
\begin{itemize}
\item Synthetic instances: \texttt{gen.py OUTDIR N SEED [medium]} (seeds: practice 1, test 2,
medium 3, training 4).
\item Built-in rules: \texttt{run\_builtin.py OUT TLIMIT WORKERS RULES INSTANCES}.
\item Network rules: \texttt{run\_net.py OUT TLIMIT WORKERS NETS [--diag K:MAX] INSTANCES}.
\item MARTA instances: \texttt{marta\_gen.py OUTDIR N SEED SERVICE 3} (pool seed 7, test seed 11;
services 5, 28, 3, 4).
\item Native plugin: \texttt{cc -O2 csrc/wspsolve.c -lscip}; runs via
\texttt{run\_c.py OUT 300 10 RULES INSTANCES} with \texttt{SCIP\_DEFAULTS=1}.
\item Figures and every number in the text: \texttt{scripts/figs.py}, which writes
\texttt{figs/summary.json}.
\end{itemize}

\section{What the AI did and what the human did}
\label{app:ai}
\begin{itemize}
\item \textbf{Human author:} chose the topic (WSP applied to exact OR), the first target (column
generation on XTW), the pivot to LP/MILP, the request for a real application, the port to C, and
the decision to publish; approved each step; checked the results.
\item \textbf{AI (Claude Opus~5.5 in Claude Code):} read the WSP post; designed and wrote all code;
wrote every network weight and note; launched all runs on the author's machine; read the
diagnostics; wrote the analysis, figures and this text; looked up and verified references.
\item \textbf{Not done by anyone:} no number was edited by hand; no run was discarded except the
instance that SCIP could not read and the superseded Python run of the final test, which was
stopped when the C plugin replaced it.
\end{itemize}

\section{Claims and evidence}
\label{app:claims}
\begin{center}\small
\resizebox{\linewidth}{!}{\begin{tabular}{@{}llll@{}}
\toprule
ID & Evidence & Log files (\texttt{results/}) & Key in \texttt{summary.json} \\
\midrule
F1 & Fig.~\ref{fig:cg} & \texttt{sweep1.jsonl} (wsp-cg) & \texttt{cg\_*} \\
F2 & Fig.~\ref{fig:syn}, \ref{fig:headroom} & \texttt{gen\_test}, \texttt{gen\_medium}, \texttt{step0} & \texttt{synthetic\_*}, \texttt{headroom} \\
F3 & Fig.~\ref{fig:agree}, \ref{fig:trained} & \texttt{sb\_train}, \texttt{vs\_trained\_*} & \texttt{sb\_agreement\_train}, \texttt{hand\_vs\_trained} \\
F4 & Fig.~\ref{fig:miplib} & \texttt{miplib}, \texttt{miplib2} & \texttt{miplib} \\
F5 & Fig.~\ref{fig:mtest}, Tab.~\ref{tab:marta} & \texttt{c\_ctrl}, \texttt{c\_test} & \texttt{marta\_test\_c}, \texttt{paired\_tests\_time} \\
F6 & Fig.~\ref{fig:mprac}, \ref{fig:mtest} & \texttt{c\_prac}, \texttt{c\_test} & \texttt{marta\_practice\_c} \\
\bottomrule
\end{tabular}}
\end{center}

\end{document}